\documentclass[12pt,fullpage]{article}

\usepackage{amsfonts,xr,graphicx,amsmath,amssymb}
\def\R{\mathbb{R}}
\def\N{\mathbb{N}}

\def\tilde{\widetilde}
\def\epsilon{\varepsilon}
   \def\ga{\gamma}  
 \def\Om{\Omega}   
\def\de{\delta} \def\De{\Delta}

\def\pa{\partial}

\def\ss{\subset}
\def\ee{\epsilon}
\def\beq{\begin{equation}}
\def\eeq{\end{equation}}

\def\div{\rm{div}}

\def\osc{\rm{osc}}
\def\length{\rm{length}}
\def\ds{\displaystyle}

\newcommand{\SE}{\setcounter{equation}{0} \section}
\newcommand{\baa}{\begin{array}}
\newcommand{\eaa}{\end{array}}
\newcommand{\ba}{\begin{eqnarray}}
\newcommand{\ea}{\end{eqnarray}}
\newtheorem{theo}{\bf Theorem}[section]
\newtheorem{lem}[theo]{\bf Lemma}
\newtheorem{pro}[theo]{\bf Proposition}

\newtheorem{rem}[theo]{\bf Remark}

\begin{document}

\title{\bf{Shear flows of an ideal fluid and elliptic equations in unbounded domains}\thanks{This work has been carried out in the framework of Archim\`ede Labex (ANR-11-LABX-0033) and of the A*MIDEX project (ANR-11-IDEX-0001-02), funded by the ``Investissements d'Avenir" French Government program managed by the French National Research Agency (ANR). The research leading to these results has also received funding from the European Research Council under the European Union's Seventh Framework Programme (FP/2007-2013) ERC Grant Agreement n.~321186~- ReaDi~- Reaction-Diffusion Equations, Propagation and Modelling and from the ANR NONLOCAL project (ANR-14-CE25-0013).}}
\author{Fran{\c{c}}ois Hamel and Nikolai Nadirashvili\\
\\
\footnotesize{Aix Marseille Universit\'e, CNRS, Centrale Marseille, I2M, UMR 7373, 13453 Marseille, France}}
\date{}
\maketitle

\begin{abstract}
We prove that, in a two-dimensional strip, a steady flow of an ideal incompressible fluid with no stationary point and tangential boundary conditions is a shear flow. The same conclusion holds for a bounded steady flow in a half-plane. The proofs are based on the study of the geometric properties of the streamlines of the flow and on one-dimensional symmetry results for solutions of some semilinear elliptic equations. Some related rigidity results of independent interest are also shown in $n$-dimensional slabs in any dimension~$n$.
\end{abstract}

\small{AMS 2000 Classification: 76B03; 35J61; 35B06; 35B53}

%%%%%%%%%%%%%%%%%%%%%%%%%%%%%%%%%%%%%%%%%%%%%%%
%%%%%%%%%%%%%%%%%%%%%%%%%%%%%%%%%%%%%%%%%%%%%%%

\SE{Introduction and main results}\label{intro}

Let $v(x)=(v^1(x),\dots,v^n(x)),\, x\in\overline{\Om}\ss\R^n,$ be a $C^2(\overline{\Om})$ velocity field of a steady flow of an ideal fluid, where $\Om$ is an open connected subset of $\R^n$. The vector field $v$ is a solution of the system of Euler equations: 
\begin{equation}\label{1}
\left\{\begin{array}{ll}
v\cdot\nabla\,v +\nabla p =0 & \mbox{in $\Om$},\vspace{3pt}\\
{\div}\ v=0 & \mbox{in $\Om$}.
\end{array} \right.
\end{equation}
The flow $v$ is called a shear flow if there is an orthogonal coordinate system $(x_1,\dots,x_n)$ in $\R^n$ such that $v$ is independent of $x_1$ and $v^2=\dots =v^n=0$. It is easy to see that the flow $v$ is a shear flow if and only if the pressure $p$ is a constant.

In this paper we give a characterization of shear flows in dimension $n=2$, in the strip $\Om_2 \ss\R^2$ defined by
\beq\label{defomega}
\Om_2 =\R\times(0,1)=\big\{ x=(x_1,x_2)\in \R^2,\ 0<x_2<1\big\}
\eeq
and in the half-plane
$$\R^2_+=\R\times(0,+\infty)=\big\{ x=(x_1,x_2)\in \R^2,\ x_2>0\big\}.$$
Throughout the paper, we denote $|\cdot|$ the Euclidean norm in $\R^m$.

Our first main result is concerned with the strip $\Omega_2$.

\begin{theo}\label{th1}
Suppose that the flow $v$ solving~\eqref{1} is defined in the closed strip $\overline{\Om_2}$ with $v^2=0$ on $\pa \Om_2$ and $\inf_{\Om_2}|v|>0$. Then $v$ is a shear flow, that is,
\beq\label{shear}
v(x)=(v^1(x_2),0)\ \hbox{ in }\overline{\Om_2}.
\eeq
\end{theo}

The condition $v^2=0$ on $\partial\Om_2$ simply means that $v$ is assumed to be tangential on the boundary of $\Om_2$.

The assumption $\inf_{\Om_2}|v|>0$ means that the flow $v$ has no stationary point in $\overline{\Om_2}$ nor limiting stationary point at infinity as $x_1\to\pm\infty$. In other words, Theorem~\ref{th1} means that any $C^2(\overline{\Om_2})$ non-shear flow which is tangential on $\partial\Om_2$ must have a stationary point in $\overline{\Om_2}$ or at infinity. These stationary points may well be in $\overline{\Om_2}$ or only at infinity. For instance, on the one hand, for any $\alpha\neq 0$, the non-shear cellular flow of the type
\beq\label{defv}
v(x)=\nabla^{\perp}\big(\sin(\alpha x_1)\sin(\pi x_2)\big)=\big(\!-\pi\sin(\alpha x_1)\cos(\pi x_2),\alpha\cos(\alpha x_1)\sin(\pi x_2)\big),
\eeq
which solves~\eqref{1} with $p(x)=(\pi^2/4)\cos(2\alpha x_1)+(\alpha^2/4)\cos(2\pi x_2)$ and is tangential on $\partial\Om_2$, has stationary points in $\overline{\Om_2}$. On the other hand, the flow
$$v(x)=\nabla^{\perp}\big(\sin(\pi x_2)\,e^{x_1}\big)=\big(\!-\pi\cos(\pi x_2)\,e^{x_1},\sin(\pi x_2)\,e^{x_1}\big),$$
which solves~\eqref{1} with $p(x)=-(\pi^2/2)e^{2x_1}$ and is tangential on $\partial\Om_2$, has no stationary point in $\overline{\Om_2}$ ($|v|>0$ in $\overline{\Om_2}$), but $\inf_{\Om_2}|v|=0$.

Lastly, we point out that, the sufficient condition $\inf_{\Om_2}|v|>0$ is obviously not equivalent to being a shear flow, since any shear flow $v(x)=(v^1(x_2),0)$ for which $v^1$ does not have a constant strict sign does not satisfy the condition $\inf_{\Om_2}|v|>0$ (however, under the conditions of Theorem~\ref{th1}, the first component $v^1$ in~\eqref{shear} has a constant strict sign in $\overline{\Om_2}$).

\begin{rem}{\rm In the assumptions of Theorem~\ref{th1}, the flow $v$ is not assumed to be a priori bounded in $\overline{\Om_2}$. However, since $v$ is $($at least$)$ continuous in $\overline{\Om_2}$ and the interval $[0,1]$ is bounded, the conclusion of Theorem~\ref{th1} implies that $v$ is necessarily bounded.}
\end{rem}

The second main result deals with the case of the half-plane $\R^2_+$.

\begin{theo}\label{th2}
Suppose that the flow $v$ solving~\eqref{1} is defined in the closed half-plane $\overline{\R^2_+}$ with $v^2=0$ on $\pa\R^2_+$ and
\beq\label{hyp}
0<\inf_{\R^2_+}|v|\le\sup_{\R^2_+}|v|<+\infty.
\eeq
Then $v$ is a shear flow, which means here that $v(x)=(v^1(x_2),0)$ in $\overline{\R^2_+}=\R\times[0,+\infty)$. In particular, the first component $v^1$ has a constant strict sign in $\overline{\R^2_+}$.
\end{theo}

Some comments are in order on the condition~\eqref{hyp}. None of the strict inequalities can be dropped for the conclusion to hold in general. More precisely, first, any cellular flow $v$ of the type~\eqref{defv}, which solves~\eqref{1} and is tangential on $\partial\R^2_+$, satisfies $\sup_{\R^2_+}|v|<+\infty$ and $\inf_{\R^2_+}|v|=0$, and it is not a shear flow. Second, the flow
$$v(x)=\nabla^{\perp}\big(x_2\cosh(x_1))=(-\cosh(x_1),x_2\sinh(x_1)),$$
which solves~\eqref{1} with $p(x)=-\cosh(2x_1)/4+x_2^2/2$ and is tangential on $\partial\R^2_+$, satisfies $\inf_{\R^2_+}|v|>0$ and $\sup_{\R^2_+}|v|=+\infty$, and it is not a shear flow. 

An interesting question would be to see whether the conclusion of Theorem~\ref{th2} still holds if the condition~\ref{hyp} is replaced by the following weaker one
$$\forall\,A>0,\ \ 0<\inf_{\R\times(0,A)}|v|\le\sup_{\R\times(0,A)}|v|<+\infty.$$
We leave it as an open problem and refer to Remark~\ref{remA} below and the end of the proof of Theorem~\ref{th2} in Section~\ref{sec3} for further comments.

The proofs of Theorems~\ref{th1} and~\ref{th2} rely on the study of the geometric properties of the streamlines of the flow $v$ and of the orthogonal trajectories of the gradient flow defined by the potential $u$ of the flow $v$ (see definition~\eqref{defu} below). The main point is to show that all streamlines of $v$ go from $-\infty$ to $+\infty$ in the direction $x_1$ (and are bounded in the direction $x_2$ in the case of the half-plane $\R^2_+$). To do so, we use a continuation argument. Therefore, the streamlines of $v$ are shown to foliate the domain and, since the vorticity $\frac{\partial v^2}{\partial x_1}-\frac{\partial v^1}{\partial x_2}$ is constant along the streamlines of the flow $v$, the potential function $u$ will be proved to satisfy a semilinear elliptic equation of the type $\Delta u=f(u)$.

To conclude the proof of Theorem~\ref{th2} (the case of the half-plane), we use some monotonicity and one-dimensional symmetry results for solutions with bounded gradient of such semilinear elliptic equations~\cite{bcn,fv}. Regarding the case of the strip $\Om_2$ (with bounded cross section), we will reduce Theorem~\ref{th1} to a new Liouville type result for these elliptic equations. More precisely, we will show the following theorem.

\begin{theo}\label{th3}
Let $f:\R\to\R$ be a locally Lipschitz continuous function and let $u$ be a $C^2(\overline{\Om_2})$ bounded solution of the equation 
\beq\label{eq}
\De u+f(u)=0
\eeq
in the strip $\Om_2$ defined in~\eqref{defomega}. Suppose that $u$ is equal to some constants on the boundary lines of $\Om_2$: $u=0$ on $\{x_2=0\}$ and $u=c$ on $\{x_2=1\}$, with $c>0$, and that
\beq\label{c12}
0<u<c
\eeq
in $\Om_2$. Then $u$ is a function of $x_2$ only, that is, $u(x_1,x_2)=\tilde{u}(x_2)$ in $\overline{\Om_2}$, and $\tilde{u}'(x_2)>0$ for all $0<x_2<1$.
\end{theo}

\begin{rem}{\rm The conclusion is sharp in the sense that $\tilde{u}'(0)$ and/or $\tilde{u}'(1)$ may well be equal to $0$ in general. For instance, the function $u(x_1,x_2)=\tilde{u}(x_2):=1-\cos(\pi x_2)$ solves~\eqref{eq} in $\Om_2$ with $f(s)=\pi^2s$, $u=0$ on $\{x_2=0\}$, $u=2$ on $\{x_2=1\}$, $0<u<2$ in $\Om_2$, $\tilde{u}'(x_2)>0$ in $(0,1)$, but $\tilde{u}'(0)=\tilde{u}'(1)=0$.}
\end{rem}

As a matter of fact, given the assumptions of Theorem~\ref{th3}, it follows from~\cite[Theorem~1.1']{bcn} applied to $u$ (resp. $c-u(x_1,1-x_2)$) that
$$u_{x_2}:=\frac{\partial u}{\partial x_2}>0\ \hbox{ in }\Big\{(x_1,x_2)\in\R^2,\ 0<x_2<\frac{1}{2}\Big\}$$
(resp. $u_{x_2}>0$ in $\big\{(x_1,x_2)\in\R^2,\ 1/2<x_2<1\big\}$). Therefore, $u_{x_2}\ge0$ in $\overline{\Om_2}$ (we point out that this monotonicity is only known in dimension $n=2$). The new result in Theorem~\ref{th3} is the fact that the monotonicity property $u_{x_2}\ge0$ in $\overline{\Om_2}$ implies that $u$ is one-dimensional, that is, $u$ is a function of $x_2$ only. Actually, it turns out that this last implication holds in any dimension $n\ge2$, as the following theorem shows:

\begin{theo}\label{th4}
Let $n\ge 2$, let $f:\R\to\R$ be a locally Lipschitz continuous function and let $u$ be a $C^2(\overline{\Om_n})$ bounded solution of the equation~\eqref{eq} in the $n$-dimensional slab $\Om_n$ defined by
\beq\label{omegan}
\Om_n =\R^{n-1}\times(0,1)=\big\{ x=(x_1,\dots,x_n)\in \R^n,\ 0<x_n<1\big\}.
\eeq
Suppose that $u$ is equal to some constants on the boundary hyperplanes of $\Om_n$: $u=0$ on $\{x_n=0\}$ and $u=c$ on $\{x_n=1\}$, with $c>0$. Suppose also that $u$ is non-decreasing with respect to the variable $x_n$, namely $u_{x_n}\ge0$ in $\overline{\Om_n}$. Then $u$ is a function of $x_n$ only, that is, $u(x_1,\dots,x_n)=\tilde{u}(x_n)$ in $\overline{\Om_n}$, and $\tilde{u}'(x_n)>0$ for all $0<x_n<1$.
\end{theo}

\begin{rem}\label{remstrict}{\rm Notice that the monotonicity assumption $u_{x_n}\ge0$ in $\overline{\Om_n}$ and the strong maximum principle imply that, for any $k\in\N$ with $k\ge2$, either
\beq\label{smp}
u(x',x_n+1/k)=u(x',x_n)\ \hbox{ for all }(x',x_n)\in\R^{n-1}\times[0,1-1/k],
\eeq
or
\beq\label{smp2}
u(x',x_n+1/k)>u(x',x_n)\ \hbox{ for all }(x',x_n)\in\R^{n-1}\times(0,1-1/k),
\eeq
where $x'=(x_1,\dots,x_{n-1})$. Indeed, the $C^2(\R^{N-1}\times[0,1-1/k])$ function $w$ defined by $w(x',x_n)=u(x',x_n+1/k)-u(x',x_n)$ is a nonnegative solution of an elliptic equation of the type $\De w+a(x)\,w=0$ in $\R^{n-1}\times(0,1-1/k)$ for some function $a\in L^{\infty}(\R^{n-1}\times(0,1-1/k))$. The case~\eqref{smp} would imply that $0=u(x',0)=u(x',1/k)=\cdots=u(x',1)=c$, a contradiction with the assumption $c>0$. Therefore,~\eqref{smp2} holds for all $k\in\N$ with $k\ge2$. Hence $u$ is actually (strictly) increasing with respect to the variable $x_n$ and
\beq\label{strict}
0<u<c\ \hbox{ in }\Om_n.
\eeq
In other words, the monotonicity assumption $u_{x_n}\ge0$ in $\overline{\Om_n}$ and the boundary conditions $u=0$ on $\{x_n=0\}$ and $u=c$ on $\{x_n=1\}$ with $c>0$ make the assumption~\eqref{c12} in $\Om_n$  redundant in the statement of Theorem~\ref{th4}. On the other hand, in Theorem~\ref{th3}, this assumption~\eqref{c12}, with strict inequalities, is essential to get the monotonicity property $u_{x_2}\ge0$ in $\Om_2$, as shown in~\cite{bcn,fs}.}
\end{rem}

In Theorem~\ref{th4}, the monotonicity condition $u_{x_n}\ge0$ in $\overline{\Om_n}$ is a sufficient condition for the one-dimensional symmetry to hold. Other sufficient conditions can be given, as the following result shows.

\begin{pro}\label{pro1}
Let $n\ge 2$, let $f:\R\to\R$ be a locally Lipschitz continuous function and let $u$ be a $C^2(\overline{\Om_n})$ bounded solution of the equation~\eqref{eq} in the $n$-dimensional slab $\Om_n$ defined in~\eqref{omegan}. Suppose that $u$ is equal to some constants on the boundary hyperplanes of $\Om_n$: $u=0$ on $\{x_n=0\}$ and $u=c$ on $\{x_n=1\}$ with $c>0$, and that condition~\eqref{c12} holds. If
\beq\label{hyp3}
\forall\,0<x_n<1,\ \ 0<\inf_{x'\in\R^{n-1}}u(x',x_n)\le\sup_{x'\in\R^{n-1}}u(x',x_n)<c
\eeq
or if
\beq\label{hyp4}
f(c)\le 0\le f(0)
\eeq
then $u$ is a function of $x_n$ only, that is, $u(x_1,\dots,x_n)=\tilde{u}(x_n)$ in $\overline{\Om_n}$, and $\tilde{u}'(x_n)>0$ for all $0<x_n<1$.
\end{pro}

In Proposition~\ref{pro1}, if condition~\eqref{hyp4} is assumed, then it follows from~\cite[Theorem~1.1]{bcn} applied to $u$ (resp. $c-u(x',1-x_n)$) that $u_{x_n}>0$ in $\big\{x\in\R^n,\ 0<x_n<1/2\big\}$ (resp. $u_{x_n}>0$ in $\big\{x\in\R^n,\ 1/2<x_n<1\big\}$). Hence $u_{x_n}\ge0$ in $\overline{\Om_n}$ and the conclusion of Proposition~\ref{pro1} in this case follows from Theorem~\ref{th4}. We point out that Theorem~\ref{th4} holds in any dimension $n\ge 2$ and without any sign assumption on $f(0)$ and $f(c)$.

The remaining sections are organized as follows. In Section~\ref{sec2}, we show that Theorem~\ref{th1} can be reduced to Theorem~\ref{th3}, that is, Theorem~\ref{th3} implies Theorem~\ref{th1}. Section~\ref{sec3} is devoted to the proof of Theorem~\ref{th2}. Regarding Theorem~\ref{th3}, as already explained, it can be deduced from earlier results of the literature and from Theorem~\ref{th4}, the latter being proved in Section~\ref{sec4} with the use of the sliding method. The proof of Proposition~\ref{pro1} is also done in Section~\ref{sec4}. From the previous paragraph, only the case of condition~\eqref{hyp3} will be considered in the proof of Proposition~\ref{pro1}.

%%%%%%%%%%%%%%%%%%%%%%%%%%%%%%%%%%%%%%%%%%%%%%%
%%%%%%%%%%%%%%%%%%%%%%%%%%%%%%%%%%%%%%%%%%%%%%%

\SE{Proof of Theorem~\ref{th1}, assuming Theorem~\ref{th3}}\label{sec2}

In this section, $v:\overline{\Om_2}\to\R^2$ is a $C^2(\overline{\Om_2})$ flow solving~\eqref{1} in the strip $\Om_2$ defined in~\eqref{defomega}. We assume that there is $\ee_0>0$ such that
$$|v(x)|\ge\ee_0>0\ \hbox{ for all }x\in\overline{\Om_2}$$
and that $v$ is tangential on the boundary, that is, $v^2=0$ on $\partial\Om_2$. Our goal is to show that $v$ is a shear flow.

Let us introduce a few important notations and definitions. First, let $u$ be a potential function of the flow $v$, that is, $u:\overline{\Om_2}\to\R$ is a $C^3(\overline{\Om_2})$ function such that
\beq\label{defu}
u_{x_1}=v^2\ \hbox{ and }\ u_{x_2}=-v^1
\eeq
in $\overline{\Om_2}$. Since $\Om_2$ is a simply connected domain, it follows that the potential function $u$ is well and uniquely defined in $\overline{\Om_2}$ up to a constant. We fix a unique function $u$ such that $u(0,0)=0$. Hence,
\beq\label{u0}
u(x_1,0)=0\ \hbox{ for all }x_1\in\R
\eeq
since $v^2=0$ on $\partial\Om_2$. For the same reason, there is a constant $c\in\R$ such that
\beq\label{u1}
u(x_1,1)=c\ \hbox{ for all }x_1\in\R.
\eeq
We will show that $u$ is bounded and ranges in $\Om_2$ between its two constant values $0$ and $c$ on $\partial\Om_2$. To do so, we need to establish some further properties of the level curves of $u$. The level curves of the function $u$ are understood as the connected components of the level sets of $u$. They are actually non-parametrized streamlines of the flow $v$ (trajectories of the flow $v$). For any $z\in\overline{\Om_2}$, we denote by $\Gamma_z$ the streamline of the flow $v$ going through $z$. In other words, $\Gamma_z$ is the level curve of $u$ containing $z$, that is the connected component of the level set $\big\{x\in\overline{\Om_2},\ u(x)=u(z)\big\}$ containing $z$. Notice that, since $|v|>0$ in $\overline{\Om_2}$ and $v^2=0$ on $\partial\Om_2$, a given streamline $\Gamma_z$ of $v$ cannot have an endpoint in $\overline{\Om_2}$ and it always admits a $C^1$ parametrization $\gamma:\R\to\R^2$ ($\gamma(\R)=\Gamma_z$) such that $|\dot{\gamma}(t)|>0$ for all $t\in\R$. Regarding the solution $\gamma_z$ of $\dot\gamma_z(t)=v(\gamma_z(t))$ with $\gamma_z(0)=z$, it is defined on a maximal open interval $(a_-,a_+)\ss\R$ which may or may not be unbounded, depending on $v$.

Second, let us consider the gradient flow
$$\dot y=\nabla u=(v^2,-v^1)$$
in $\overline{\Om_2}$, whose trajectories are orthogonal to the streamlines of the flow $v$. Since $|\nabla u |=|v|\geq \ee_0>0$ in $\overline{\Om_2}$, the following lemma holds immediately.

\begin{lem}\label{lem1}
Let $G$ be a nonempty open set included in $\Om_2$. Assume that $u$ is bounded in $G$ and let $g:\R\to\R^2$ be a parametrization of a trajectory $\Sigma$ of the gradient flow $\dot y=\nabla u$ in $G$. Then $\Sigma$ has a finite length such that
$$
{\length}(\Sigma) \leq \frac{{\osc}_G\,u}{\ee_0}:=\frac{\sup_Gu-\inf_Gu}{\ee_0}
$$
and it terminates on the boundary $\pa G$ of $G$ $($that is, ${\rm{dist}}(g(t),\partial G)=\inf_{x\in\partial G}|g(t)-x|\to0$ as $t\to\pm\infty$$)$.
\end{lem}

\begin{rem}\label{rembounded}{\rm 
Notice that Lemma~\ref{lem1} and the assumption $|\nabla u |=|v|\geq \ee_0>0$ imply that, for any $z\in G$, the solution $\sigma_z$ of $\dot\sigma_z(t)=\nabla u(\sigma_z(t))$ with $\sigma_z(0)=0$ and $\sigma_z(t)\in G$ is defined on a maximal interval $(a_-,a_+)\ss\R$ which is necessarily bounded (however, the curve $\sigma_z((a_-,a_+))$ always admits parametrizations defined in $\R$).}
\end{rem}

The following three lemmas deal with some continuity (with respect to $z$) and asymptotic properties of the streamlines $\Gamma_z$ of the flow $v$. In what follows, for any $r>0$ and $z\in\overline{\Om_2}$, we set
$$B_{z,r}=\big\{ x\in\overline{\Om_2},\ |x-z|<r\big\}.$$

\begin{lem}\label{lem2}
Let $z_1\in\overline{\Om_2}$ be given. For any $\ee>0$, there is $\de>0$ such that, for any $z_2\in B_{z_1,\de}$ and $z_3\in B_{z_1,\de}$, $\Gamma_{z_3}$ lies in an $\ee$-neighborhood of  $\Gamma_{z_2}$.
\end{lem}

\noindent{\bf{Proof.}} We fix $z_1\in\overline{\Om_2}$. For any $r>0$, we define
\beq\label{defGz1r}
G_{z_1,r}=\mathop{\bigcup}_{z\in B_{z_1,r}}\Gamma_z.
\eeq
By definition of $u$ and the streamlines $\Gamma_z$, there holds
$${\osc}_{G_{z_1,r}}u={\osc}_{B_{z_1,r}}u.$$ 
Hence, there exists a constant $C>0$ only depending on $z_1$ and $u$ such that, for all $0<r<1$,
$$
{\osc}_{G_{z_1,r}}u\leq C\,r.
$$\par
Since $|\nabla u(z_1)|=|v(z_1)|>0$, there is $r_1\in(0,1]$ such that, for every $r\in(0,r_1)$, every level set of $u$ has only one connected component in $B_{z_1,r}$ and can be written as the graph of a function depending on the variable $x\cdot v(z_1)$. Hence, for every $r\in(0,r_1)$, $z_2\in B_{z_1,r}$ and $z_3\in B_{z_1,r}\backslash\Gamma_{z_2}$, there holds $u(z_3)\neq u(z_2)$. Moreover, for any such $r$, $z_2$ and $z_3$, let us define
\beq\label{defG'}
G'_{z_2,z_3,r}=\mathop{\bigcup}_{z\in B_{z_1,r},\ \min(u(z_2),u(z_3))<u(z)<\max(u(z_2),u(z_3))}\Gamma_z,
\eeq
which is a nonempty open set included in $G_{z_1,r}\cap\Om_2$, with $\partial G'_{z_2,z_3,r}=\Gamma_{z_2}\cup\Gamma_{z_3}$. Notice that these properties hold even if $B_{z_1,r}$ intersects $\partial\Om_2$: in that case, since $\nabla u$ is orthogonal to $\partial\Omega_2$ on $\partial\Omega_2$, the values $0$ on $\{x_2=0\}$ or $c$ on $\{x_2=1\}$ would be a global minimum or maximum of $u$ in $B_{z_1,r}$, still for $r>0$ small enough.\par
Now, for any $r\in(0,r_1)$, any $z_2\in B_{z_1,r}$, any $z_3\in B_{z_1,r}\backslash \Gamma_{z_2}$ and any $z_4\in \Gamma_{z_3}$, let $\Sigma_{z_4}\ss G'_{z_2,z_3,r}$ be the trajectory of the gradient flow $\dot y=\nabla u$ with an end point at $z_4$ (this gradient curve $\Sigma_{z_4}$ is well defined in $G'_{z_2,z_3,r}$ since by definition $\nabla u(z_4)\neq(0,0)$ is orthogonal to $\Gamma_{z_3}$ at $z_4$). By Lemma~\ref{lem1}, the second end point of $\Sigma_{z_4}$ lies on $\partial G'_{z_2,z_3,r}$ and hence on $\Gamma_{z_2}$ (it can not lie on $\Gamma_{z_3}$ since $u$ is strictly monotone along $\Sigma_{z_4}$). Moreover,
$${\length}(\Sigma_{z_4})\leq\frac{{\osc}_{G'_{z_2,z_3,r}}u}{\ee_0}\le\frac{{\osc}_{G_{z_1,r}}u}{\ee_0}\le\frac{C\,r}{\ee_0}.$$
Hence, for every $r\in(0,r_1)$, every $z_2\in B_{z_1,r}$ and every $z_3\in B_{z_1,r}\backslash \Gamma_{z_2}$, the distance between any point $z_4\in \Gamma_{z_3}$ and  $\Gamma_{z_2}$ is less than $Cr/\ee_0$ and the lemma follows.\hfill$\Box$

\begin{lem}\label{unbounded}
For any given streamline $\Gamma$ of the flow $v$ and any $C^1$ parametrization $\gamma=(\gamma^1,\gamma^2):\R\to\Gamma$ such that $|\dot{\gamma}(t)|>0$ for all $t\in\R$, there holds $|\ga^1(t)|\to+\infty $ as $|t|\to+\infty$. 
\end{lem}

\noindent{\bf{Proof.}} Assume that the conclusion does not hold. Then there are $y\in\overline{\Om_2}$ and a sequence $(\tau_n)_{n\in\N}$ such that $|\tau_n|\to+\infty$ and $\gamma(\tau_n)\to y$ as $n\to+\infty$. Since $|v(y)|>0$, it follows that, for each $k\in\N$ large enough, there is $t_k\in\R$ such that
\beq\label{tildetk}
v(y)\cdot(\gamma(t_k)-y)=0,\ \ |\gamma(t_k)-y|\le k^{-1}\ \hbox{ and }\ |t_k|>|t_{k-1}|.
\eeq\par
If there were some integers $k$ and $l$ such that $t_k<t_l$ and $\gamma(t_k)=\gamma(t_l)$, then the open set $\omega$ surrounded by $\gamma([t_k,t_l])$ would be nonempty (since $|\dot{\gamma}(t)|>0$ for all $t\in\R$). By definition of $\gamma$, the function $u$ is constant on the curve $\gamma([t_k,t_l])=\partial\omega$. Thus, $u$ has either an interior minimum or an interior maximum in $\omega$, which is ruled out since $u$ has no critical point. Therefore, the points $\gamma(t_k)$ are pairwise distinct.\par
Now, from~\eqref{tildetk}, for each $k\in\N$ large enough, the nonzero vector $\gamma(t_k)-\gamma(t_{k-1})$ is parallel to the nonzero vector $\nabla u(y)$ (which is orthogonal to $v(y)$), whence
$$\frac{\gamma(t_k)-\gamma(t_{k-1})}{|\gamma(t_k)-\gamma(t_{k-1})|}=\pm\frac{\nabla u(y)}{|\nabla u(y)|}.$$
One infers that, for each $k$ large enough, there is $\theta_k\in[0,1]$ such that
$$\baa{rcl}
\ds\frac{|u(\gamma(t_k))-u(\gamma(t_{k-1}))|}{|\gamma(t_k)-\gamma(t_{k-1})|} & = & \ds\frac{|\nabla u\big(\theta_k\gamma(t_k)+(1-\theta_k)\gamma(t_{k-1})\big)\cdot(\gamma(t_k)-\gamma(t_{k-1}))|}{|\gamma(t_k)-\gamma(t_{k-1})|}\vspace{3pt}\\
& = & \ds\frac{|\nabla u\big(\theta_k\gamma(t_k)+(1-\theta_k)\gamma(t_{k-1})\big)\cdot\nabla u(y)|}{|\nabla u(y)|}.\eaa$$
But the left-hand side is equal to $0$ since $u(\gamma(t_k))=u(\gamma(t_{k-1}))$ by definition of $\gamma$, whereas the right-hand side converges to $|\nabla u(y)|>0$ as $k\to+\infty$, a contradiction. The proof of Lemma~\ref{unbounded} is thereby complete.\hfill$\Box$\break

Now, we say that a streamline $\Gamma$ of the flow $v$ (that is, a level curve of $u$) is {\it regular} if it has a parametrization $\gamma:\R\to\Gamma$ such that
\beq\label{regular}
\ga^1(t)\to\pm\infty\ \hbox{ as }t\to\pm\infty,
\eeq
with $\ga(t)=(\ga^1(t),\ga^2(t))$.

\begin{lem}\label{lem3}
All streamlines of the flow $v$ are regular.
\end{lem}

\noindent{\bf{Proof.}} Let $a$ be an arbitrary real number. Since $v^2=0$ on the line $\Gamma=\{(x_1,0),\ x_1\in\R\}$ and $v$ is continuous with $|v|\ge\ee_0>0$ in $\overline{\Om_2}$, the line $\Gamma$ is equal to the streamline $\Gamma_{(a,0)}$ and it is regular. Let us consider
$$E=\big\{b\in[0,1]:\Gamma_{(a,b)}\hbox{ is regular}\big\}.$$
The set $E$ is not empty, since $0\in E$ (notice also that $1\in E$ for the same reason as for $0$). We have to show that $E=[0,1]$. To do so, let us prove that $E$ is both closed and relatively open in $[0,1]$.\par
First, let $b$ be in $E$, that is, the streamline $\Gamma_{(a,b)}$ is regular. By Lemma~\ref{lem2}, there is $\de>0$ such that $\Gamma_{(a,b)}$ lies in a $1$-neighborhood of $\Gamma_{(a,b')}$ for every $b'\in[b-\delta,b+\delta]\cap[0,1]$. For any such $b'$, by Lemma~\ref{unbounded}, the streamline $\Gamma_{(a,b')}$ has a $C^1$ parametrization $\gamma=(\gamma^1,\gamma^2):\R\to\Gamma_{(a,b')}$ such that $|\dot\gamma(t)|>0$ for every $t\in\R$ and $|\gamma^1(t)|\to+\infty$ as $|t|\to+\infty$. In particular, up to changing $t$ into $-t$, one can assume without loss of generality that
$$\gamma^1(t)\to-\infty\ \hbox{ as }t\to-\infty.$$
Since the streamline $\Gamma_{(a,b)}$ is assumed to be regular and since it lies in a $1$-neighborhood of $\Gamma_{(a,b')}$, it follows immediately that $\gamma^1(t)\to+\infty$ as $t\to+\infty$. Therefore, $\Gamma_{(a,b')}$ is regular for every $b'\in[b-\delta,b+\delta]\cap[0,1]$, and $E$ is relatively open in $[0,1]$.\par
Second, let $b\in\overline{E}$. By Lemma~\ref{lem2}, there is $\de>0$ such that $\Gamma_{(a,b')}$ lies in a $1$-neighborhood of $\Gamma_{(a,b)}$ for every $b'\in[b-\delta,b+\delta]\cap[0,1]$. Since $b\in\overline{E}$, there is $b'\in E\cap[b-\delta,b+\delta]$. In particular, the streamline $\Gamma_{(a,b')}$ is regular. Thus, as in the previous paragraph, the streamline $\Gamma_{(a,b)}$ is regular too. In other words, $b\in E$ and $E$ is closed.\par
As a conclusion, $E$ is equal to the whole interval $[0,1]$. Since $a$ was an arbitrary real number, one infers that all streamlines $\Gamma_z$ (for any $z\in\overline{\Om_2}$) of the flow $v$ are regular.\hfill$\Box$\break

We recall that $u$ is constant along the lines $\{x_2=0\}$ and $\{x_2=1\}$, and satisfies~\eqref{u0} and~\eqref{u1}, after normalization. Furthermore, since $v^2=0$ on $\partial\Om_2$ and $|v|>0$ in $\overline{\Om_2}$, one can assume in the sequel without loss of generality that
$$v^1(0,0)<0,$$
even if it means changing $v$ into $-v$. 

\begin{lem}\label{lem4}
The function $u$ is bounded in $\overline{\Om_2}$, $c>0$ and there holds
$$0<u<c\ \hbox{ in }\Om_2.$$
\end{lem}

\noindent{\bf{Proof.}} 
Let $M\in[0,+\infty)$ be defined as
$$M=\max_S|u|,$$
where $S$ is the segment $S=\{(0,x_2),\ x_2\in[0,1]\}$. By Lemma~\ref{lem3} all level curves of $u$ are regular and then intersect $S$. Therefore, $|u|$ is bounded by $M$ in $\overline{\Om_2}$.\par
Let now $\sigma$ be the solution of
\beq\label{defsigma}
\dot\sigma(t)=\nabla u(\sigma(t)),
\eeq
taking values in $\overline{\Om_2}$, with
\beq\label{sigma00}
\sigma(0)=(0,0).
\eeq
Since $u$ is bounded in $\Om_2$ and $|\nabla u|\ge\ee_0>0$ in $\Om_2$ with $\nabla u(0,0)=(v^2(0,0),-v^1(0,0))=(0,-v^1(0,0))$ and $-v^1(0,0)>0$, it follows from Lemma~\ref{lem1} and Remark~\ref{rembounded} that the function $\sigma$ is defined on a maximal time interval $[0,\tau]$ with $\tau\in(0,+\infty)$ and that
\beq\label{defSigma}
\Sigma:=\sigma([0,\tau])
\eeq
has finite length and ends on $\partial\Om_2$. Since the function
\beq\label{deftheta}
\theta:t\mapsto u(\sigma(t))
\eeq
is (strictly) increasing in the interval $[0,\tau]$ ($\dot\theta(t)=|\nabla u(\sigma(t))|^2\ge\ee_0^2>0$ for all $t\in[0,\tau]$), the second end point $\sigma(\tau)$ of $\Sigma$ lies on $\{x_2=1\}$, that is,
$$\sigma(\tau)=(\xi,1)$$
for some $\xi\in\R$. Furthermore, since $\sigma(0)=u(0,0)=0$ and $\sigma(\tau)=u(\xi,1)=c$, one infers that $c>0$ and that $0<u(\sigma(t))<c$ for all $t\in(0,\tau)$.\par
Let now $x$ be an arbitrary point in $\Om_2$. By Lemma~\ref{lem3}, the streamline $\Gamma_x$ is regular and therefore it intersects $\Sigma=\sigma([0,\tau]\big)$. In other words, there is $t_x\in[0,\tau]$ such that $\sigma(t_x)\in\Gamma_x$, hence
$$u(\sigma(t_x))=u(x)$$
by definition of $u$ and $\Gamma_x$. Notice that $t_x$ is unique since $t\mapsto u(\sigma(t))$ is (strictly) increasing on $[0,\tau]$. Furthermore, the streamline $\Gamma_x$ lies entirely in the open set $\Om_2$ (it can not intersect $\partial\Om_2=\{x_2=0\}\cup\{x_2=1\}$ since both lines $\{x_2=0\}$ and $\{x_2=1\}$ are themselves streamlines of $v$). Thus, $\sigma(t_x)\in\Om_2$. Hence $t_x\in(0,\tau)$ and $0<u(\sigma(t_x))<c$, that is,
$$0<u(x)<c.$$
Since $x\in\Om_2$ was arbitrary, the proof of Lemma~\ref{lem4} is thereby complete.\hfill$\Box$\break

\noindent{\bf{End of the proof of Theorem~\ref{th1}.}} Let $\sigma:[0,\tau]\to\overline{\Om_2}$ be the solution of~\eqref{defsigma} with $\sigma(0)=(0,0)$ and $\sigma(\tau)=(\xi,1)$ for some $\xi\in\R$. Since $u$ is of class $C^3(\overline{\Om_2})$, the function $\sigma$ is of class $C^3([0,\tau])$. The function $\theta:[0,\tau]\to[0,c]$ defined by~\eqref{deftheta} satisfies $\dot\theta(t)=|\nabla u(\sigma(t))|^2>0$ in $[0,\tau]$ and is therefore a $C^3$ diffeomorphism from $[0,\tau]$ to $[0,c]$. Let us now define
$$\baa{rcl}
f:[0,c] & \!\!\to\!\! & \R\vspace{3pt}\\
s & \!\!\mapsto\!\! & f(s):=-\Delta u\big(\sigma(\theta^{-1}(s))\big),\eaa$$
that is,
$$f(\theta(t))=-\Delta u(\sigma(t))$$
for all $t\in[0,\tau]$. The function $f$ is of class $C^1$ in $[0,c]$ (by extending it to $f(0)$ in $(-\infty,0)$ and by $f(c)$ in $(c,+\infty)$, it then becomes a Lipschitz continuous function defined in $\R$).\par
Let us finally show that $u$ is a classical solution of the elliptic equation $\Delta u+f(u)=0$ in $\overline{\Om_2}$. To do so, observe first that $\Delta u=v^2_{x_1}-v^1_{x_2}$ in $\overline{\Om_2}$, and that the $C^1(\overline{\Om_2})$ function $v^2_{x_1}-v^1_{x_2}$ satisfies
$$v\cdot\nabla(v^2_{x_1}-v^1_{x_2})=0\ \hbox{ in }\overline{\Om_2},$$
by~\eqref{1}. Therefore, the function $v^2_{x_1}-v^1_{x_2}$ is constant along any streamline of $v$, that is, along any level curve of $u$. Let now $x$ denote any arbitrary point in $\overline{\Om_2}$. As in the proof of Lemma~\ref{lem4}, the regular streamline $\Gamma_x$ intersects $\Sigma=\sigma\big([0,\tau]\big)$ and there is a unique $t_x\in[0,\tau]$ such that $\sigma(t_x)\in\Gamma_x$ and
$$\theta(t_x)=u(\sigma(t_x))=u(x).$$
Finally, since the function $v^2_{x_1}-v^1_{x_2}$ is constant on the streamline $\Gamma_x$ containing both $x$ and $\sigma(t_x)$, one infers from the definitions of $\theta$ and $f$ that
\beq\label{equ}
\Delta u(x)=v^2_{x_1}(x)-v^1_{x_2}(x)=v^2_{x_1}(\sigma(t_x))-v^1_{x_2}(\sigma(t_x))=\Delta u(\sigma(t_x))=-f(\theta(t_x))=-f(u(x)).
\eeq\par
As a conclusion, the function $u$ is a classical solution of~\eqref{eq} in $\overline{\Om_2}$ with $u=0$ on $\{x_2=0\}$, $u=c$ on $\{x_2=1\}$ and $0<u<c$ in $\Om_2$, for some Lipschitz continuous function $f:\R\to\R$. Theorem~\ref{th3} implies that $u$ depends only on the variable $x_2$, that is $u(x)=\tilde{u}(x_2)$ in $\overline{\Om_2}$, for some $C^3([0,1])$ function $\tilde{u}$ such that $\tilde{u}'(x_2)>0$ for all $x_2\in(0,1)$. In other words,
$$v(x)=(-\tilde{u}'(x_2),0)\ \hbox{ in }\overline{\Om_2}.$$
Notice finally that $v^1(x)=-\tilde{u}'(x_2)$ has a constant sign in $\overline{\Om_2}$, including the boundary $\partial\Om_2$, since $|v|$ is continuous and does not vanish in $\overline{\Om_2}$. The proof of Theorem~\ref{th1} is thereby complete.\hfill$\Box$

\begin{rem}{\rm It follows from the conclusion of Theorem~\ref{th1} that, for any $z=(z_1,z_2)\in\overline{\Omega_2}$, the streamline $\Gamma_z$ is nothing but the line $\{(x_1,x_2)\in\R^2,\ x_2=z_2\}$. Furthermore, the gradient curve $\Sigma$ defined in~\eqref{defSigma} is equal to the segment $S=\{(0,x_2),\ x_2\in[0,1]\}$.}
\end{rem}

%%%%%%%%%%%%%%%%%%%%%%%%%%%%%%%%%%%%%%%%%%%%%%%
%%%%%%%%%%%%%%%%%%%%%%%%%%%%%%%%%%%%%%%%%%%%%%%

\SE{Proof of Theorem~\ref{th2}}\label{sec3}

The proof of Theorem~\ref{th2} follows the same scheme as that of Theorem~\ref{th1}, apart from some additional observations induced by the unboundedness of $\R^2_+$ in the direction $x_2$. For the sake of clarity, we preferred to put the two proofs in two different sections.\par
In this section, $v:\overline{\R^2_+}\to\R^2$ is a $C^2(\overline{\R^2_+})$ flow solving~\eqref{1} in the half-plane $\R^2_+$. We assume that there are $0<\ee_0\le M<+\infty$ such that
$$0<\ee_0\le|v(x)|\le M\ \hbox{ for all }x\in\overline{\R^2_+}$$
and that $v$ is tangential on the boundary, that is, $v^2=0$ on $\partial\R^2_+=\{(x_1,0),\ x_1\in\R\}$. Our goal is to show that $v$ is a shear flow.

As in Section~\ref{sec2}, let $u$ be a potential function of the flow $v$, that is, $u:\overline{\R^2_+}\to\R$ is a $C^3(\overline{\R^2_+})$ function satisfying~\eqref{defu} in $\overline{\R^2_+}$. The function $u$ is uniquely fixed by the normalization $u(0,0)=0$, hence
$$u(x_1,0)=0\ \hbox{ for all }x_1\in\R,$$
since $v^2=0$ on $\partial\R^2_+$. The level curves of $u$ are the streamlines $\Gamma_z$ of the flow $v$. Since $|v|>0$ in $\overline{\R^2_+}$ and $v^2=0$ on $\partial\R^2_+$, any given streamline $\Gamma_z$ of $v$ can not have an endpoint in $\overline{\R^2_+}$ and it has a $C^1$ parametrization $\gamma:\R\to\R^2$ ($\gamma(\R)=\Gamma_z$) such that $|\dot{\gamma}(t)|>0$ for all $t\in\R$. Since $|\nabla u|=|v|\geq \ee_0>0$ in $\R^2_+$, Lemma~\ref{lem1} holds immediately, with $\R^2_+$ in place of $\Omega_2$.

For any $r>0$ and $z\in\overline{\R^2_+}$, we define the restricted ball $B_{z,r}$ as
$$B_{z,r}=\big\{ x\in\overline{\R^2_+},\ |x-z|<r\big\}.$$
It is obvious to see that Lemma~\ref{lem2} holds, with $\R^2_+$ in place of $\Omega_2$.

Now, for any $z\in\overline{\R^2_+}$, we say that the streamline $\Gamma_z$ is {\it vertically bounded} if there is a real number $A_z>0$ such that
\beq\label{vertical}
\Gamma_z\,\subset\,\R\times[0,A_z].
\eeq
Notice that, for any $x_1\in\R$, the streamline $\Gamma_{(x_1,0)}=\partial\R^2_+$ is vertically bounded. We shall prove that all streamlines $\Gamma_z$ are vertically bounded and regular, in the sense of~\eqref{regular}. To do so, we first observe that the analogue of Lemma~\ref{unbounded} holds for vertically bounded streamlines.

\begin{lem}\label{unbounded2}
For any given vertically bounded streamline $\Gamma$ of the flow $v$ and any $C^1$ parametrization $\gamma=(\gamma^1,\gamma^2):\R\to\Gamma$ such that $|\dot{\gamma}(t)|>0$ for all $t\in\R$, there holds $|\ga^1(t)|\to+\infty $ as $|t|\to+\infty$. 
\end{lem}

\noindent{\bf{Proof.}} It is identical to that of Lemma~\ref{unbounded} up to replacing $\overline{\Omega_2}$ by $\R\times[0,A]$, where $A>0$ is such that $\Gamma\subset\R\times[0,A]$.\hfill$\Box$\break

Lemma~\ref{lem3} can now be extended as follows.

\begin{lem}\label{lem5}
All streamlines $\Gamma_z$ of the flow $v$ are regular in the sense of~\eqref{regular} and vertically bounded in the sense of~\eqref{vertical}.
\end{lem}

\noindent{\bf{Proof.}} Let $a$ be an arbitrary real number. Since $v^2=0$ on the line $\Gamma=\{(x_1,0),\ x_1\in\R\}=\partial\R^2_+$ and $v$ is continuous with $|v|\ge\ee_0>0$ in $\overline{\R^2_+}$, $\Gamma=\Gamma_{(a,0)}$ and this streamline is regular and vertically bounded. Let us consider
\beq\label{defE}
E=\big\{b\in[0,+\infty),\ \Gamma_{(a,b)}\hbox{ is regular and vertically bounded}\big\}.
\eeq
The set $E$ is not empty, since $0\in E$. Let us now show that $E$ is both closed and relatively open in $[0,+\infty)$.\par
First, let $b$ be in $E$, that is, the streamline $\Gamma_{(a,b)}$ is regular and vertically bounded. By Lemma~\ref{lem2} (with $\R^2_+$ in place of $\Om_2$), there is $\de>0$ such that $\Gamma_{(a,b)}$ lies in a $1$-neighborhood of $\Gamma_{(a,b')}$ and $\Gamma_{(a,b')}$ lies in a $1$-neighborhood of $\Gamma_{(a,b)}$ for every $b'\in[b-\delta,b+\delta]\cap[0,+\infty)$. For any such $b'$, the streamline $\Gamma_{(a,b')}$ is vertically bounded, since so is $\Gamma_{(a,b)}$. By Lemma~\ref{unbounded2}, $\Gamma_{(a,b')}$ has a $C^1$ parametrization $\gamma=(\gamma^1,\gamma^2):\R\to\Gamma_{(a,b')}$ such that $|\dot\gamma(t)|>0$ for every $t\in\R$ and $|\gamma^1(t)|\to+\infty$ as $|t|\to+\infty$. Since $\Gamma_{(a,b)}$ is regular, one then concludes as in the proof of Lemma~\ref{lem3} that $\Gamma_{(a,b')}$ is regular. Therefore, $b'\in E$. Hence, $[b-\delta,b+\delta]\cap[0,+\infty)\,\ss\,E$ and $E$ is relatively open in $[0,+\infty)$.\par
Second, let $b\in\overline{E}$. By Lemma~\ref{lem2} (with $\R^2_+$ in place of $\Om_2$), there is $\de>0$ such that $\Gamma_{(a,b)}$ lies in a $1$-neighborhood of $\Gamma_{(a,b')}$ and $\Gamma_{(a,b')}$ lies in a $1$-neighborhood of $\Gamma_{(a,b)}$ for every $b'\in[b-\delta,b+\delta]\cap[0,+\infty)$. Since $b\in\overline{E}$, there is $b'\in E\cap[b-\delta,b+\delta]$. In particular, the streamline $\Gamma_{(a,b')}$ is vertically bounded and regular. Thus, as in the previous paragraph, the streamline $\Gamma_{(a,b)}$ is vertically bounded and regular too. In other words, $b\in E$ and $E$ is closed.\par
As a conclusion, $E$ is equal to the whole interval $[0,+\infty)$. Since $a$ was an arbitrary real number, the conclusion of Lemma~\ref{lem5} follows.\hfill$\Box$\break

We recall that $u$ is constant along the lines $\Gamma=\{(x_1,0),\ x_1\in\R\}=\partial\R^2_+$, with $u=0$ on $\partial\R^2_+$ after normalization. As in Section~\ref{sec2}, since $v^2=0$ on $\partial\R^2_+$ and $|v|>0$ in $\overline{\R^2_+}$, one can assume in the sequel without loss of generality that
$$v^1(0,0)<0.$$

\begin{lem}\label{lem6}
There holds
$$u>0\hbox{ in }\R^2_+.$$
\end{lem}

\noindent{\bf{Proof.}} Let $\sigma$ be the solution of~\eqref{defsigma} and~\eqref{sigma00} taking values in $\overline{\R^2_+}$, that is $\dot\sigma(t)=\nabla u(\sigma(t)),
$ and $\sigma(0)=(0,0)$. Since $\nabla u(0,0)=(v^2(0,0),-v^1(0,0))=(0,-v^1(0,0))$ and $-v^1(0,0)>0$, the function $\sigma$ is defined on a maximal interval $I$ of the type $[0,\tau]$ with $\tau\in(0,+\infty)$, or $[0,\tau)$ with $\tau\in(0,+\infty]$. The gradient curve $\Sigma=\sigma(I)$ cannot end at an (interior) point in $\R^2_+$, since $|\nabla u|=|v|>0$ in $\R^2_+$. It cannot end on $\partial\R^2_+$ either, since the function
\beq\label{deftheta2}
\theta:t\mapsto u(\sigma(t))
\eeq
is (strictly) increasing in $I$ ($\dot\theta(t)=|\nabla u(\sigma(t))|^2\ge\ee_0^2>0$ for all $t\in I$) and $u$ is constant on $\partial\R^2_+$. Therefore, $\Sigma$ has infinite length and, since $|\dot\sigma(t)|=|\nabla u(\sigma(t))|\le M$ for all $t\in I$, one infers that $I$ is infinite, that is $I=[0,+\infty)$ and
\beq\label{defSigma2}
\Sigma=\sigma([0,+\infty)).
\eeq
Moreover, $u(\sigma(t))>u(\sigma(0))=u(0,0)=0$ for all $t\in(0,+\infty)$.\par
Write $\sigma(t)=(\sigma^1(t),\sigma^2(t))$ for all $t\in[0,+\infty)$. We claim that
\beq\label{claim2}
\limsup_{t\to+\infty}\sigma^2(t)=+\infty.
\eeq
Indeed, otherwise, the function $t\mapsto\sigma^2(t)$ would be bounded in $[0,+\infty)$ and there would exist $A>0$ such that $\Sigma\subset\R\times[0,A]$. But $u_{x_2}=-v^1$ is bounded in $\R\times[0,A]$ (since $v$ is actually assumed to be bounded in $\R^2_+$) and $u=0$ on $\{x_2=0\}$. Thus, $u$ is bounded in $\R\times[0,A]$, hence $\theta(t)=u(\sigma(t))$ is bounded in $[0,+\infty)$. But $\dot\theta(t)=|\nabla u(\sigma(t))|^2\ge\ee_0^2>0$ for all $t\in[0,+\infty)$, leading to a contradiction. As a consequence, the claim~\eqref{claim2} is proved.\par
Let now $x$ be an arbitrary point in $\R^2_+$. By Lemma~\ref{lem5}, the streamline $\Gamma_x$ is regular and vertically bounded. Therefore, by~\eqref{claim2} and $\sigma(0)=(0,0)$, $\Gamma_x$ intersects $\Sigma=\sigma([0,+\infty)\big)$. In other words, there is $t_x\in[0,+\infty)$ such that $\sigma(t_x)\in\Gamma_x$, hence
$$u(\sigma(t_x))=u(x).$$
The real number $t_x$ is unique since $t\mapsto u(\sigma(t))$ is (strictly) increasing. Furthermore, the streamline $\Gamma_x$ lies entirely in the open set $\R^2_+$ (it can not intersect $\partial\R^2_+=\{x_2=0\}$ since the line $\{x_2=0\}$ is itself a streamline of $v$). Thus, $\sigma(t_x)\in\R^2_+$. Hence $t_x\in(0,+\infty)$ and $u(\sigma(t_x))>0$, that is, $u(x)>0$. The proof of Lemma~\ref{lem6} is thereby complete.\hfill$\Box$

\begin{rem}\label{remA}{\rm It is straightforward to see that, once Lemma~\ref{lem5} is proved, the conclusion of Lemma~\ref{lem6} holds if, instead of~\eqref{hyp}, one only assumes that
\beq\label{hypA}
\forall\,A>0,\ \ 0<\inf_{\R\times(0,A)}|v|\le\sup_{\R\times(0,A)}|v|<+\infty,
\eeq
even if the interval $[0,\tau)$ might be bounded in that case. Lemma~\ref{lem1} and~\ref{lem2} could be extended under this hypothesis, provided that $G\subset\R^2_+$ in Lemma~\ref{lem1} is vertically bounded and $\Gamma_{z_1}\subset\overline{\R^2_+}$ in Lemma~\ref{lem2} is vertically bounded. However, it is unclear to prove Lemma~\ref{lem5} with~\eqref{hypA} instead of~\eqref{hyp} (the proof of the closedness of the set $E$ defined in~\eqref{defE} is not clear).}
\end{rem}

\noindent{\bf{End of the proof of Theorem~\ref{th2}.}} Let $\sigma:[0,+\infty)\to\overline{\R^2_+}$ be the solution of~\eqref{defsigma} and~\eqref{sigma00}. Since $u$ is of class $C^3(\overline{\R^2_+})$, the function $\sigma$ is of class $C^3([0,+\infty))$. The function $\theta:[0,+\infty)\to[0,+\infty)$ defined by~\eqref{deftheta2} satisfies $\dot\theta(t)=|\nabla u(\sigma(t))|^2\ge\ee_0^2>0$ in $[0,+\infty)$ with $\theta(0)=u(\sigma(0))=u(0,0)=0$ and is therefore a $C^3$ diffeomorphism from $[0,+\infty)$ to itself. Let us now define
$$\baa{rcl}
f:[0,+\infty) & \!\!\to\!\! & \R\vspace{3pt}\\
s & \!\!\mapsto\!\! & f(s):=-\Delta u\big(\sigma(\theta^{-1}(s))\big),\eaa$$
that is, $f(\theta(t))=-\Delta u(\sigma(t))$ for all $t\in[0,+\infty)$. The function $f$ is of class $C^1$ in $[0,+\infty)$ (by extending it to $f(0)$ in $(-\infty,0)$, it is then a locally Lipschitz continuous function defined in $\R$).\par
Let us finally show that $u$ is a classical solution of the elliptic equation $\Delta u+f(u)=0$ in $\overline{\R^2_+}$. As in the proof of Theorem~\ref{th1}, $\Delta u=v^2_{x_1}-v^1_{x_2}$ and $v\cdot\nabla(v^2_{x_1}-v^1_{x_2})=0$ in $\overline{\R^2_+}$. Therefore, the function $v^2_{x_1}-v^1_{x_2}$ is constant along any streamline of $v$. Let now $x$ denote an arbitrary point in $\overline{\R^2_+}$. From Lemma~\ref{lem5} and~\eqref{claim2}, the regular and vertically bounded streamline $\Gamma_x$ intersects $\Sigma=\sigma\big([0,+\infty)\big)$ and there is a unique $t_x\in[0,+\infty)$ such that $\sigma(t_x)\in\Gamma_x$ and $\theta(t_x)=u(\sigma(t_x))=u(x)$. Finally, since $v^2_{x_1}-v^1_{x_2}$ is constant on the streamline $\Gamma_x$ containing both $x$ and $\sigma(t_x)$, one infers as in~\eqref{equ} that
$$\Delta u(x)=-f(u(x)).$$\par
As a conclusion, the function $u$ is a classical solution of~\eqref{eq} in $\overline{\R^2_+}$ with $u=0$ on $\partial\R^2_+=\{x_2=0\}$ and $u>0$ in $\R^2_+$, for some locally Lipschitz-continuous function $f:\R\to\R$. For any $A>0$, the function $u$ is positive and bounded in the strip $\R\times(0,A)$, it is then a classical solution of an equation of the type $\Delta u+f_A(u)=0$ in $\R\times[0,A]$ for some globally Lipschitz continuous function $f_A:[0,+\infty)\to\R$ such that $f_A=f$ on the image of $\R\times[0,A]$ by $u$. It then follows from~\cite[Theorem~1.1']{bcn} that
$$u_{x_2}>0\hbox{ in }\R\times\Big(0,\frac{A}{2}\Big),$$
hence $u_{x_2}>0$ in $\R^2_+$ since $A>0$ can be arbitrary (see also~\cite[Theorem~2]{d1}, ~\cite[Remark~2]{d2} and~\cite[Theorem~1.1]{fs} for related results). Finally, since $|\nabla u|=|v|$ is bounded in $\R^2_+$ (the second strict inequality in~\eqref{hyp} is used here), one concludes from~\cite[Theorem~1.3]{fv} (see also~\cite[Theorems~1.6]{fs}) that $u$ depends only on the variable $x_2$, that is
$$u(x)=\tilde{u}(x_2)\hbox{ in }\overline{\R^2_+},$$
for some $C^3([0,+\infty))$ function $\tilde{u}$ such that $\tilde{u}'(x_2)>0$ for all $x_2\in(0,+\infty)$. In other words,
$$v(x)=(-\tilde{u}'(x_2),0)\ \hbox{ in }\overline{\R^2_+}.$$
Furthermore, $v^1(x)=-\tilde{u}'(x_2)$ has a constant sign in $\overline{\R^2_+}$, including the boundary $\partial\R^2_+$, since $|v|$ is continuous and does not vanish in $\overline{\R^2_+}$. The proof of Theorem~\ref{th2} is thereby complete.\hfill$\Box$

\begin{rem}{\rm It follows from the conclusion of Theorem~\ref{th2} that, for any $z=(z_1,z_2)\in\overline{\R^2_+}$, the streamline $\Gamma_z$ is nothing but the line $\{(x_1,x_2)\in\R^2,\ x_2=z_2\}$. Furthermore, the gradient curve $\Sigma$ defined in~\eqref{defSigma2} is equal to the half-line $\{(0,x_2),\ x_2\in[0,+\infty)\}$.}
\end{rem}

%%%%%%%%%%%%%%%%%%%%%%%%%%%%%%%%%%%%%%%%%%%%%%%
%%%%%%%%%%%%%%%%%%%%%%%%%%%%%%%%%%%%%%%%%%%%%%%

\SE{Proof of Theorem~\ref{th4} and Proposition~\ref{pro1}}\label{sec4}

This section is devoted to the proof of the one-dimensional symmetry results concerned with the elliptic equation~\eqref{eq}. As already emphasized in Section~\ref{intro}, Theorem~\ref{th3} follows from Theorem~\ref{th4} and from earlier results in the literature. It only remains to show Theorem~\ref{th4} and Proposition~\ref{pro1}.\hfill\break

\noindent{\bf{Proof of Theorem~\ref{th4}.}} Let $\Om_n$ be the $n$-dimensional slab defined in~\eqref{omegan}, with $n\ge 2$. Let $f:\R\to\R$ be a locally Lipschitz continuous function and let $u$ be a $C^2(\overline{\Om_n})$ bounded solution of the equation~\eqref{eq} in $\Omega_n$ with $u=0$ on $\{x_n=0\}$, $u=c$ on $\{x_n=1\}$ and $c>0$. We also assume that $u_{x_n}\ge0$ in $\overline{\Om_n}$.\par
In order to prove that $u$ depends on the variable $x_n$ only, we will show that it is increasing in any direction $\xi=(\xi_1,\cdots,\xi_n)=(\xi',\xi_n)$ with $\xi'=(\xi_1,\cdots,\xi_{n-1})\in\R^{n-1}$ and $\xi_n>0$. We fix such a vector $\xi$ and we will use a sliding method~\cite{bn}, by sliding $u$ in the direction $\xi$ and comparing the shifted function to the function $u$ itself. Namely, for any $\tau\in(0,1/\xi_n)$, we define
\beq\label{omegatau}
\Omega^{\tau}=\R^{n-1}\times(0,1-\tau\xi_n)
\eeq
and
\beq\label{wtau}
w^{\tau}(x)=u(x+\tau\xi)-u(x)\ \hbox{ for }x\in\overline{\Omega^{\tau}}.
\eeq
Our aim is to show that $w^{\tau}>0$ in $\overline{\Omega^{\tau}}$ for all $\tau\in(0,1/\xi_n)$.\par
First, it follows from standard elliptic estimates up to the boundary that $\nabla u$ is bounded in $\overline{\Omega_n}$. Since $u(x',0)=0<c=u(x',1)$ for all $x'\in\R^{n-1}$, one gets the existence of $\tau^*\in(0,1/\xi_n)$ such that
$$w^{\tau}>0\ \hbox{ in }\overline{\Omega^{\tau}}\ \hbox{ for all }\tau\in(\tau^*,1/\xi_n).$$
Let us now define
\beq\label{deftau*}
\tau_*=\inf\big\{\tau\in(0,1/\xi_n),\ w^{\tau'}>0\ \hbox{ in }\overline{\Omega^{\tau'}}\hbox{ for all }\tau'\in(\tau,1/\xi_n)\big\}.
\eeq
There holds $0\le\tau_*\le\tau^*<1/\xi_n$.\par
Let us assume by contradiction that $\tau_*>0$. One has $w^{\tau}>0$ in $\overline{\Omega^{\tau}}$ for all $\tau\in(\tau_*,1/\xi_n)$ and $w^{\tau_*}\ge0$ in $\overline{\Omega^{\tau_*}}$ by continuity. Moreover, there are a sequence $(\tau_k)_{k\in\N}$ in $(0,\tau_*]$ converging to $\tau_*$ and a sequence $(x^k)_{k\in\N}$ of points in $\R^n$ such that
$$x^k\in\overline{\Omega^{\tau_k}}\ \hbox{ and }\ w^{\tau_k}(x^k)\le0$$
for all $k\in\N$. Write $x^k=({x'}^k,x^k_n)$ with ${x'}^k\in\R^{n-1}$ and $x^k_n\in[0,1-\tau_k\xi_n]$, and define
\beq\label{defukwk}\left\{\baa{ll}
u_k(x)=u(x'+{x'}^k,x_n) & \hbox{for }x\in\overline{\Omega_n},\vspace{3pt}\\
w_k(x)=w^{\tau_*}(x'+{x'}^k,x_n)=u_k(x+\tau_*\xi)-u_k(x) & \hbox{for }x\in\overline{\Omega^{\tau_*}}.\eaa\right.
\eeq
From standard elliptic estimates, since $f$ is locally Lipschitz continuous and $u$ is a $C^2(\overline{\Om_n})$ bounded solution of~\eqref{eq}, the sequences of the functions $u_k$ and $w_k$ are actually bounded in $C^{2,\alpha}_{loc}(\overline{\Om_n})$ and $C^{2,\alpha}_{loc}(\overline{\Om^{\tau_*}})$ respectively, for every $\alpha\in[0,1)$. Therefore, up to extraction of a subsequence, the functions $u_k$ converge in $C^2_{loc}(\overline{\Om_n})$ to a $C^2(\overline{\Om_n})$ bounded solution $U$ of~\eqref{eq}. By passing to the limit, one gets that $U=0$ on $\{x_n=0\}$, $U=c$ on $\{x_n=1\}$ and $U_{x_n}\ge0$ in $\overline{\Omega_n}$ (hence, $0\le U\le c$ in $\overline{\Om_n}$). It follows then from Remark~\ref{remstrict}, as in~\eqref{strict}, that
\beq\label{strict2}
0<U<c\ \hbox{ in }\Omega_n.
\eeq\par
Now, by definition of $w_k$, the functions $w_k$ then converge in $C^2_{loc}(\overline{\Omega^{\tau_*}})$ to the function $W$ defined in $\overline{\Om^{\tau_*}}$ by
$$W(x)=U(x+\tau_*\xi)-U(x)\ \hbox{ for }x\in\overline{\Omega^{\tau_*}}.$$
Furthermore, $W\ge0$ in $\overline{\Om^{\tau_*}}$ (since so is $w^{\tau_*}$). Up to extraction of a subsequence, one can assume that $x^k_n\to\tilde{x}_n$ as $k\to+\infty$ for some $\tilde{x}_n\in[0,1-\tau_*\xi_n]$. Let us set $\tilde{x}=(0,\cdots,0,\tilde{x}_n)\in\overline{\Omega^{\tau_*}}$. Lastly, $w^{\tau_k}(x^k)\le0$ means that
$$u({x'}^k+\tau_k\xi',x^k_n+\tau_k\xi_n)\le u({x'}^k,x^k_n),$$
that is, $u_k(\tau_k\xi',x^k_n+\tau_k\xi_n)\le u_k(0,x^k_n)$. By passing to the limit as $k\to+\infty$, one infers that $U(\tau_*\xi',\tilde{x}_n+\tau_*\xi_n)\le U(0,\tilde{x}_n)$, that is $U(\tilde{x}+\tau_*\xi)\le U(\tilde{x})$, i.e. $W(\tilde{x})\le0$. But $W\ge0$ in $\overline{\Omega^{\tau_*}}$. Hence, $W(\tilde{x})=0$, that is
$$U(\tilde{x}+\tau_*\xi)=U(\tilde{x}).$$\par
Several cases can occur, whether $\tilde{x}$ be on $\partial\Omega^{\tau_*}=\{x_n=0\}\cup\{x_n=1-\tau_*\xi_n\}$ or in $\Omega^{\tau_*}$. On the one hand, if $\tilde{x}_n=0$, then $U(\tilde{x})=0$. Hence, $U(\tilde{x}+\tau_*\xi)=0$. But $\tilde{x}+\tau_*\xi\in\Omega_n$ (it is an interior point) since $0<\tilde{x}_n+\tau_*\xi_n=\tau_*\xi_n<1$, contradicting~\eqref{strict2}. On the other hand, if $\tilde{x}_n=1-\tau_*\xi_n$, then $U(\tilde{x}+\tau_*\xi)=c$. Hence, $U(\tilde{x})=c$. But $\tilde{x}\in\Omega_n$ since $0<\tilde{x}_n=1-\tau_*\xi_n<1$, contradicting again~\eqref{strict2}. Therefore, $\tilde{x}\in\Omega^{\tau_*}$ (it is an interior point). But the function $W$ is a nonnegative solution of an equation of the type $\Delta W+c(x)\,W=0$ in $\Omega^{\tau_*}$, for some function $c\in L^{\infty}(\Omega^{\tau_*})$. The strong maximum principle implies that $W=0$ in $\Omega^{\tau_*}$, and then on $\partial\Omega^{\tau_*}$ by continuity. This leads to a contradiction as in the beginning of this paragraph.\par
Finally, the assumption $\tau_*>0$ is ruled out. In other words, $\tau_*=0$ and  
$$u(x+\tau\xi)>u(x)\ \hbox{ for all }x\in\overline{\Omega^{\tau}}\hbox{ and for all }0<\tau<1/\xi_n.$$
This means that $u$ is increasing in any direction $\xi=(\xi',\xi_n)\in\R^{n-1}\times\R$ such that $\xi_n>0$. By continuity, for any $\xi'\in\R^{n-1}$, $u$ then is nondecreasing in the direction $\xi=(\xi',0)$. So is it in the direction $\xi=(-\xi',0)$. As a consequence, $u$ does not depend on the direction $(\xi',0)$ and, since $\xi'\in\R^{n-1}$ is arbitrary, $u$ depends only on the variable $x_n$, that is $u(x)=\tilde{u}(x_n)$ for some $C^2([0,1])$ function $\tilde{u}$. The function $\tilde{u}$ is actually increasing in $[0,1]$ from Remark~\ref{remstrict}. Lastly, for any $h\in(0,1/2]$, the $C^2([0,h])$ function $z$ defined by
$$z(x_n)=\tilde{u}(2h-x_n)-\tilde{u}(x_n)$$
is nonnegative and it satisfies an equation of the type $z''(x_n)+d(x_n)\,z(x_n)=0$ in $[0,h]$ for some bounded function $d$. Furthermore, $z(h)=0$, and $z>0$ in $[0,h)$ since $\tilde{u}$ is increasing in $[0,1]$. Hopf lemma implies that $z'(h)<0$, that is, $\tilde{u}'(h)>0$ for any $h\in(0,1/2]$. Similarly, by working with the function $c-\tilde{u}(1-x_n)$, one infers that $\tilde{u}'(h)>0$ for all $h\in[1/2,1)$. The proof of Theorem~\ref{th4} is thereby complete.\hfill$\Box$\break

\noindent{\bf{Proof of Proposition~\ref{pro1}.}} Let $n\ge 2$, let $f:\R\to\R$ be a locally Lipschitz continuous function and let $u$ be a $C^2(\overline{\Om_n})$ bounded solution of~\eqref{eq} in $\Omega_n$ with $u=0$ on $\{x_n=0\}$, $u=c$ on $\{x_n=1\}$ and $c>0$. As emphasized in Section~\ref{intro}, we only need to consider the case where $u$ satisfies~\eqref{hyp3}. As in the proof of Theorem~\ref{th4}, we will show that $u$ is increasing in any direction $\xi=(\xi',\xi_n)\in\R^n_+$ with $\xi'\in\R^{n-1}$ and $\xi_n>0$.\footnote{We could also prove with the same method that $u$ is increasing in $x_n$ and then use the conclusion of Theorem~\ref{th4}. But we preferred to prove directly the monotonicity in the direction $\xi$, since the notations will be the same as in the proof of Theorem~\ref{th4}.}\par
We fix such a vector $\xi$. For any $\tau\in(0,1/\xi_n)$, we define $\Omega^{\tau}$ and $w^{\tau}$ as in~\eqref{omegatau} and~\eqref{wtau}. Since $u(x',0)=0<c=u(x',1)$ for all $x'\in\R^{n-1}$, one gets as in the proof of Theorem~\ref{th4} the existence of $\tau^*\in(0,1/\xi_n)$ such that $w^{\tau}>0$ in $\overline{\Omega^{\tau}}$ for all $\tau\in(\tau^*,1/\xi_n)$. We then define $\tau_*\in[0,\tau^*]$ as in~\eqref{deftau*} and assume by contradiction that $\tau_*>0$. One has $w^{\tau}>0$ in $\overline{\Omega^{\tau}}$ for all $\tau\in(\tau_*,1/\xi_n)$ and $w^{\tau_*}\ge0$ in $\overline{\Omega^{\tau_*}}$. Moreover, there are a sequence $(\tau_k)_{k\in\N}$ in $(0,\tau_*]$ converging to $\tau_*$ and a sequence $(x^k)_{k\in\N}$ of points in $\R^n$ such that $x^k\in\overline{\Omega^{\tau_k}}$ and $w^{\tau_k}(x^k)\le0$ for all $k\in\N$. Define $u_k$ and $w_k$ as in~\eqref{defukwk}. Up to extraction of a subsequence, the functions $u_k$ converge in $C^2_{loc}(\overline{\Om_n})$ to a $C^2(\overline{\Om_n})$ bounded solution $U$ of~\eqref{eq} such that $U=0$ on $\{x_n=0\}$ and $U=c$ on $\{x_n=1\}$. Furthermore,~$U$ satisfies
$$\forall\,0<x_n<1,\ \ 0<\inf_{x'\in\R^{n-1}}U(x',x_n)\le\sup_{x'\in\R^{n-1}}U(x',x_n)<c$$
since $u$ satisfies this assumption~\eqref{hyp3} and this condition is invariant by translation in the directions $x'$. In particular, $0<U<c$ in $\Omega_n$, that is, $U$ still satisfies~\eqref{strict2}. One then gets a contradiction exactly as in the proof of Theorem~\ref{th4}.\par
Therefore, $\tau_*=0$ and $u(x+\tau\xi)>u(x)$ for all $\xi=(\xi',\xi_n)$ with $\xi_n>0$ and for all $\tau\in(0,1/\xi_n)$ and $x\in\overline{\Omega^{\tau}}$. As in the proof of Theorem~\ref{th4}, one concludes that $u(x)=\tilde{u}(x_n)$ only depends on $x_n$ and that $\tilde{u}'(x_n)>0$ for all $x_n\in(0,1)$.\hfill$\Box$

%%%%%%%%%%%%%%%%%%%%%%%%%%%%%%%%%%%%%%%%%%%%%%%
%%%%%%%%%%%%%%%%%%%%%%%%%%%%%%%%%%%%%%%%%%%%%%%

\end{document}